\documentclass[12pt]{article}
\usepackage{amsfonts,amssymb,version}

\def\AA{{\mathbb A}}
\def\B{{\mathcal B}}
\def\Cc{{\mathcal C}}

\def\K{{\mathcal K}}
\def\M{{\mathcal M}}
\def\mm{\mathfrak m}

\def\OO{{\mathcal O}}
\def\PP{{\mathbb P}}
\def\pp{\mathfrak p}

\def\qq{\mathfrak q}

\def\S{{\mathcal S}}
\def\W{{\mathcal W}}
\def\Z{\mathbb Z}
\def\Zz{{\mathcal Z}}

\def\div{\mathop{\rm div}\nolimits} 
\def\id{\mathop{\rm id}\nolimits}

\def\ord{\mathop{\rm ord}\nolimits}

\def\spec{\mathop{\rm Spec}\nolimits}

\def\supp{\mathop{\rm supp}\nolimits}

\def\toto{\stackrel{\to}{\scriptstyle \to}}

\let\da\downarrow
\let\hra\hookrightarrow
\let\ov\overline

\makeatletter \def\l@section{\@dottedtocline{1}{0em}{1.2em}} \makeatother

\begin{document}

\centerline{\Large\bf On flat pullbacks for Chow groups}

\bigskip

\centerline{\bf Nitin Nitsure}

{\footnotesize Retired Professor, School of Mathematics, 
Tata Institute of Fundamental Research,
Homi Bhabha Road, Mumbai 400 005,
India. email: nitsure@gmail.com}

\pagestyle{empty}

\begin{abstract}
It is a fundamental property of the Chow groups of algebraic schemes that
they are contra-functorial with respect to flat morphisms between schemes.
While the pullback homomorphism is easy to define at the level of algebraic
cycles, the crucial step is to show that the pullback 
of cycles preserves rational equivalence, so that it descends to the Chow groups.
The purpose of this note is to give a natural sheaf theoretic proof of
the preservation of rational equivalence under flat pullback on cycles. 
\end{abstract}

\centerline{\small Key words: Algebraic cycles, flat morphisms, rational equivalence, Chow groups.}

\centerline{\small Math. Subj. Class. MSC2010 : 14C15 (Primary), 14C25, 14A15 (Secondary).}

\bigskip

{\large \bf  1. Introduction}

\medskip

If $f : X \to Y$ is a flat morphism of algebraic schemes over a base field $k$, such that 
$f$ is of constant relative dimension $d$, then it
induces a graded homomorphism $f^*: Z_*(Y) \to Z_*(X)$ of degree $d$
between the respective $\Z$-graded groups $Z_*(Y)$ and $Z_*(X)$ of algebraic cycles,
defined by sending the generator $[V]\in Z_n(Y)$ corresponding to an $n$-dimensional closed
subvariety $V\subset Y$ to the fundamental cycle $[f^{-1}V] \in Z_{n+d}(X)$ of the 
closed subscheme $f^{-1}V \subset X$, which is of pure dimension $n+d$ if non-empty.
These pullback homomorphisms can be easily seen to satisfy 
(i) $\id_X^*$ is the identity on $Z_*(X)$, and (ii) if 
$f : X \to Y$ and $g : Y \to Z$ are flat morphisms of constant relative dimensions 
then $(g\circ f)^* = f^* \circ g^*$.     
This makes $Z_*$ a contravariant functor from the category of algebraic schemes and flat morphisms
of constant relative dimensions to the category of graded abelian groups and graded homomorphisms. 
What is more difficult is to
show that if $f : X \to Y$ is a flat morphism of constant relative dimension $d$,
then $f^*: Z_*(Y) \to Z_*(X)$
maps the graded subgroup $Rat_*(Y) \subset Z_*(Y)$ formed by 
algebraic cycles on $Y$ that are rationally equivalent to $0$, 
into the corresponding graded subgroup $Rat_*(X)\subset Z_*(X)$, thereby 
inducing a graded homomorphism $f^* : A_*(Y) \to A_*(X)$ of degree $d$ between the respective graded 
Chow groups, which are the quotients $A_*(Y) = Z_*(Y)/Rat_*(Y)$ and $A_*(X) = Z_*(X)/Rat_*(X)$.
This is the flat pullback for Chow groups in the title of this note.

The proof of the preservation of rational equivalence by flat morphisms, given
by Fulton in Chapter 1 of [3], is based on three main ingredients.
The first ingredient is the fundamental theorem of Chevalley [1] that the push-forward
of algebraic cycles by proper morphisms preserves rational equivalence. 
The second ingredient is an alternative description of rational equivalence of algebraic cycles 
on a variety $V$ by representing a rational function on $V$ via a flat
morphism $\tilde{V}\to \PP^1$
on a suitable variety $\tilde{V}$ that is proper birational over $V$, 
and then taking the push-forward on $V$ of the pullback on $\tilde{V}$
of the algebraic cycle $[P_0] - [P_{\infty}]$ on $\PP^1$.  
The third ingredient is a lemma on how the push-forwards of algebraic cycles
by proper morphisms commute 
with the pullbacks of algebraic cycles by flat morphisms of constant relative dimensions
on algebraic cycles, in a cartesian square.

The purpose of this note is to give a more `natural' proof of 
the preservation of rational equivalence by flat morphisms, by using the
sheaf theoretic properties of algebraic cycles and of invertible meromorphic functions.
The sheaf of meromorphic functions $\K_X$ on a scheme $X$ is a somewhat neglected topic in most textbooks
of algebraic geometry, hence we recall some necessary basics in Section 2
(the experienced reader can skip this section).
In Section 3, we introduce the site $\Cc$ of algebraic schemes and flat morphisms
of constant relative dimensions, equipped with the big Zariski topology.
Invertible meromorphic functions, and algebraic cycles, form sheaves $\K^{\times}$ and $\Zz_*$ on $\Cc$.
In Section 4, we introduce the sheaf $\W$ of Weil divisors on the full sub site
$\Cc^{pure}$ of $\Cc$ consisting of pure dimensional schemes, and a homomorphism
$\Phi : \K^{\times} \to \W$ over $\Cc^{pure}$, using which we complete the proof of the
preservation of rational equivalence.

\pagestyle{myheadings}
\markright{Nitin Nitsure: On flat pullbacks for Chow groups}

\bigskip

{\large \bf  2. Sheaves of meromorphic functions}

\medskip

For any ring $R$, let $S(R)\subset R$ denote set of all non-zero divisors, and let $Q(R)$
denote the total quotient ring $R[S(R)^{-1}]$ of $R$. 
For any locally ringed space $(X,\OO_X)$ there is a subsheaf of sets $\S_X\subset \OO_X$
whose sections over an open subset $U\subset X$ are all $f\in \OO_X(U)$ such that for each $x\in U$ the
germ $f_x $ lies in $S(\OO_{X,x})$
(in particular, $\S_X(U) \subset S(\OO_X(U))$, but the inclusion can sometimes be proper). 
The sheaf $\K_X$ of meromorphic functions on $X$ is the sheaf of $\OO_X$-algebras
which is the localization $\K_X = \OO_X[\S_X^{-1}]$ of $\OO_X$ obtained by inverting $\S_X$
(see EGA IV.20 [2]). If $X$ is a locally noetherian scheme, then
the sections of $\S_X$ and $\K_X$ over any affine open subscheme $U=\spec A \subset X$, 
and the stalks $\S_{X,x}$ and $\K_{X,x}$ at any $x = \pp \in \spec A$,
have the following direct algebraic description (see Mumford [5] Chapter 9, Section 1): 
$\S_X(U) = S(A)$, $\S_{X,x} = S(\OO_{X,x}) = S(A_{\pp})$, $\K_X(U) = Q(A)$ and
$\K_{X,x} = Q(\OO_{X,x}) = Q(A_{\pp})$.
If $\spec B \subset \spec A$ are affine open subschemes of $X$, then
the corresponding ring homomorphism $A\to B$ is flat, so it maps $S(A)$ into $S(B)$
and induces a homomorphism $Q(A)\to Q(B)$. These are the restriction maps
from $\spec A$ to $\spec B$ for the sections of the sheaves
$\S_X$ and $\K_X$.
Let $\K_X^{\times}$ denote the sheaf of multiplicative groups formed by the 
units in $\K_X$. It satisfies $\K_X^{\times}(\spec A) = Q(A)^{\times}$
and $(\K_X^{\times})_x = \K_{X,x}^{\times} = Q(\OO_{X,x})^{\times} = Q(A_{\pp})^{\times}$
for any affine open subscheme $\spec A\subset X$ and point $x = \pp \in \spec A$.

If $f: X \to Y$ is a flat morphism of locally noetherian schemes, then
as flat homomorphisms of rings preserve non zero-divisors, the homomorphism
$f^{\sharp}: f^{-1}\OO_Y \to \OO_X$ induces homomorphisms 
$f^{\sharp}: f^{-1}\S_Y \to \S_X$, $f^{\sharp}: f^{-1}\K_Y \to \K_X$ and
$f^{\sharp}: f^{-1}\K_Y^{\times} \to \K_X^{\times}$.
If $f: X\to Y$ and $g: Y\to Z$ are both flat then 
$f^{\sharp}\circ (f^{-1}(g^{\sharp})) = (g\circ f)^{\sharp}$
on $\S$, $\K$ and $\K^{\times}$.

The above description of the sections and stalks of $\S_X$, $\K_X$ and $\K_X^{\times}$ 
does not in general extend to non locally noetherian schemes (see Kleiman [4]). 

\medskip

{\bf Support $|r|$ and restriction map $res_{X,Y}$. }
The quotient sheaf $\K_X^{\times}/\OO_X^{\times}$ is called the sheaf of Cartier divisors
on $X$. The image of any $r\in \Gamma(X,\K_X^{\times})$ in $\Gamma(X,\K_X^{\times}/\OO_X^{\times})$
is called the Cartier divisor of $r$, and is denoted by $\div(r)$. The set-topological support
of $\div(r)$, which is a closed subset of the underlying
topological space $|X|$ of $X$, will be called as the support of $r$, and will be denoted by $|r|$
(the support $\supp(r)$ of the section $r\in \Gamma(X,\K_X^{\times})$ in the
usual sheaf theoretic sense, which is by definition 
the set of all points $x$ for which the germ $r_x \ne 1 \in \K^{\times}_{X,x}$, 
is contained in $|r|$ but $\supp(r)\ne |r|$ in general).
If $U = X - |r|$ is regarded as an open subscheme of $X$, then $r|_U \in \Gamma(U, \OO_U^{\times})$.
If $X = \spec A$ where $A$ is noetherian, then we can write $r = s/t$ where $s, t \in S(A)$, and hence
$|r| \subset V(a) \cup V(b)$ where $V(a), V(b)$ are the closed subsets defined by $a=0$ and $b=0$.
If $\pp$ is a minimal prime of $A$, then $\pp \subset Zerodiv(A)$, and so $a\not\in \pp$,
$b\not\in \pp$, showing that $\pp \in U = X - |r|$. For any $r\in \Gamma(X,\K_X^{\times})$
for a locally noetherian scheme $X$, the above argument shows that $|r|$ does not contain
any height $0$ points of $X$. Hence if $Y\subset X$ is 
an irreducible component of $X$ regarded as a closed integral subscheme, then
$U\cap Y \ne \emptyset$, and therefore the regular
function $r|_U\in \Gamma(U, \OO_U^{\times})$ further restricts to a section
$r_Y \in \Gamma(Y\cap U, \OO_Y^{\times})$. As $Y$ is integral,
$\Gamma(Y\cap U, \OO_Y^{\times}) \subset K(Y)^{\times}$ where $K(Y)$ denotes the function field of
$Y$, which in turn equals $\K_Y(Y)$. Hence we have a restriction homomorphism
$res_{X,Y} : \K^{\times}_X(X) \to \K^{\times}_Y(Y) : r \mapsto r_Y$ for any irreducible component $Y$ of $X$.  
If $U\subset X$ is an open subscheme such that $Y\cap U \ne \emptyset$, then 
$Y\cap U$ is an irreducible component of $U$, and it follows from the above definition of 
restrictions that $r_Y|_{Y\cap U} = (r|_U)_{Y\cap U} \in K_{Y\cap U}(Y\cap U)$, an equality
that will be useful later. We will denote $r_Y|_{Y\cap U}$ and $(r|_U)_{Y\cap U}$ by
the common symbol $r|_{Y \cap U}$.

\medskip

{\footnotesize For completeness, we include the following algebraic lemma 
which is used in proving the above statements about sections and stalks of $\S_X$ and $\K_X$.

{\bf Lemma 1.} {\it 
(1) For any ring $A$, let the subset $\S(A) \subset A$ consist of all $a\in A$ such that
for every prime ideal $\pp \subset A$, the image of $a$ under the localization
homomorphism $A\to A_{\pp}$ is a non zero-divisor in $A_{\pp}$.
Then $\S(A) \subset S(A)$. If the ring $A$ is noetherian, then $\S(A) = S(A)$. 
     
(2) If $A$ is noetherian and $I \subset A$ an ideal such that
$Ann(I)=0$ then $S(A) \cap I \ne \emptyset$,
where $Ann(I)\subset A$ denotes the annihilator ideal of $I$.
    
(3) Let $f_1,\ldots, f_n$ generate the unit ideal in $A$. 
If $A$ is noetherian then the diagrams
$$S(A) \to {\textstyle \prod}_i\, S(A_{f_i}) \toto {\textstyle \prod}_{(j,k)} \, S(A_{f_jf_k})
\mbox{ and } 
Q(A) \to {\textstyle \prod}_i\, Q(A_{f_i}) \toto {\textstyle \prod}_{(j,k)} \, Q(A_{f_jf_k}),$$ 
obtained by applying the functors
$S$ and $Q$ to the exact diagram
$A \to  {\textstyle \prod}_i\, A_{f_i} \toto {\textstyle \prod}_{(j,k)} A_{f_jf_k }$
of rings and flat homomorphisms, are exact.
} 

\medskip

{\bf Proof.} (1) The inclusion $\S(A) \subset S(A)$ is obvious. 
If $A$ is noetherian, its set of associated prime ideals $Ass(A)$ 
is finite and the set $Zerodiv(A)$ of all zero divisors in $A$ equals $\cup\, Ass(A)$. 
If $\qq$ is any prime ideal in $A$, then
$Ass(A_{\qq}) = \{\pp A_{\qq}\,|\, \pp \in Ass(A) \mbox{ and } \pp \subset \qq \}$.
Hence if $a \in S(A) = A - \cup\, Ass(A)$, then $a/1 \in S(A_{\qq})$,
showing $S(A) \subset \S(A)$. 

(2) Suppose to the contrary that $I \subset Zerodiv(A) = \cup Ass(A)$.
As $Ass(A)$ is a finite set of prime ideals, by prime avoidance there exists $\pp \in Ass(A)$
such that $I\subset \pp$. Being an associated prime, $\pp$ equals $Ann(x)$ the annihilator
ideal of some nonzero element $x\in A$. Hence $x\in Ann(I)$ showing that $Ann(I) \ne 0$.

(3) A proof of this, assuming (2), is given on the page 61 of Mumford [5].  
\hfill$\square$

} 

\bigskip

\bigskip

\bigskip

{\large \bf  3. Sheaves of algebraic cycles}

\medskip

We will follow the notation and terminology
used in Fulton's {\it Intersection Theory} [3],
and we will assume that the reader is familiar with the material covered in the Chapter 1 of that book.
From here onwards `schemes' will mean algebraic schemes over a chosen base field $k$,
and all morphisms will be assumed to be over $k$. 
The phrase `flat morphism' will mean a flat morphism $f : X \to Y$ of schemes 
which is of a constant relative dimension. A `variety' will mean an integral scheme, 
not assumed to be separated. If $V$ is a variety, $R(V)$ will denote the
field of rational functions over $V$. In terms of the sheaf of meromorphic functions on $V$,
note that $R(V) = \K_V(V)$, and $R(V)^{\times} = \K_V^{\times}(V)$.

Let $\Cc_{k, flat}$ denote the category whose objects are all algebraic $k$-scheme, and morphisms
are all flat $k$-morphisms of constant relative dimensions. We give it the big Zariski
topology, in other words, an open cover of an object $X$ is a Zariski open cover of $X$.  
This defines a site $\Cc_{k, flat, Zar}$ which we will simply denote by $\Cc$.
We define a sheaf $\K^{\times}$ of abelian groups on this site as follows.
For any object $X$ of $\Cc$ we define $\K^{\times}(X) = \K_X^{\times}(X)$,
and for any morphism $f: X \to Y$, we define the restriction homomorphism
$\K^{\times}(Y) \to \K^{\times}(X)$ to be the group homomorphism
induced by $f^{\sharp}: f^{-1}\K_Y \to \K_X$ (it matters here that $f$ is flat,
but it does not matter here whether $f$ has a constant relative dimension or not). 

We next define a presheaf $\Zz_*$ of $\Z$-graded abelian groups and graded homomorphisms
on $\Cc$ as follows.
For any object $X$, we define $\Zz_*(X) = Z_*(X)$, and for any morphism
$f: X \to Y$ in $\Cc$, which is by definition a flat $k$-morphism of a constant relative dimension
say $d$, we define the restriction homomorphism $\Zz_*(Y)\to \Zz_*(X)$
to be the graded homomorphism of degree $d$
induced by the flat pullbacks $f^* : Z_n(Y) \to Z_{n+d}(X)$ under $f$. 
If $V \subset X$ is a closed subvariety of dimension $n$, then $f^{-1}(V)$ is a closed
subscheme of $X$ of pure dimension $n+d$, and by definition 
$f^*[V]$ is the fundamental cycle
$[f^{-1}(V)] \in Z_{n+d}(X)$. Here, recall that if $X$ is a scheme and $W\subset X$ is a closed subscheme
of pure dimension $m$, then the fundamental cycle $[W] \in Z_m(W)\subset Z_m(X)$ is defined by
$$[W] = {\textstyle \sum}_{W_i}\, \ell(\OO_{W,W_i}) [W_i]$$ 
where the sum is indexed by all irreducible components $W_i$ of $W$, and 
$\OO_{W,W_i}$ denotes the stalk of $\OO_W$ at the generic point of $W_i$, which is an artin local ring,
with length denoted by $\ell(\OO_{W,W_i})$. The map $f^*$ so defined on the generators $[V]$ extends
uniquely to the free abelian group $Z_n(Y)$
to define a homomorphism $f^* : Z_n(Y) \to Z_{n+d}(X)$. By [3] Lemma 1.7.1,
for any closed subscheme $W$ of $Y$ of pure dimension $n$, we have
$f^*[W] = [f^{-1}W]$. From this, it follows that for any two arrows $f: X \to Y$ and $g: Y\to Z$
in $\Cc$, we have $(g\circ f)^* = f^* \circ g^*$, and it is clear that $\id_X^* = \id_{Z^*(X)}$. 
Hence $\Zz_*$ is a presheaf on $\Cc$.

\bigskip

{\bf Proposition 2.} {\it The presheaf $\Zz_*$ is a sheaf on the site $\Cc$.}

\medskip

{\bf Proof.} Let $\alpha = \sum_i m_i [V_i] \in Z_n(X)$ where $V_i$ are closed subvarieties of
$X$ of dimension $n$, with $V_i \ne V_j$ for $i\ne j$, and with $m_i \ne 0$ for each $i$.
Let $(U_{\lambda})$ be a Zariski open cover
of $X$. If $\alpha|_{U_{\lambda}} = 0$ for each $\lambda$, then we will show that $\alpha =0$.
If $V_i \cap U_{\lambda} \ne \emptyset$, then $V_i$ is uniquely recovered from $V_i \cap U_{\lambda}$
as its closure. Hence as the $V_i$ are distinct, if 
$V_i \cap U_{\lambda}  = V_j \cap U_{\lambda}$ for $i\ne j$ then we must have
$V_i \cap U_{\lambda}  = V_j \cap U_{\lambda} = \emptyset$. 
Hence there can be no cancellations amongst the coefficients $m_i$, and so if $\alpha|_{U_{\lambda}} = 0$
then each $V_i \cap U_{\lambda}  = \emptyset$. If $\alpha|_{U_{\lambda}} = 0$ for each $\lambda$,
then as $(U_{\lambda})$ is a cover of $X$, it follows that $\alpha = 0$ (the empty sum). 
This shows that $\Zz_*$ is a separated presheaf on $\Cc$.

As $\Zz_*$ is a separated presheaf and as any object $X$ of $\Cc$ is quasicompact,
to verify the gluing condition for $\Zz_*$, it is enough to consider a union
$X = U' \cup U''$ of two open subschemes, together with $n$-cycles $\alpha' = 
\sum_{i=1}^p m_i' [V_i'] \in Z_n(U')$ and 
$\alpha'' = \sum_{j=1}^q m_j'' [V_j''] \in Z_n(U'')$ where $p,q\ge 0$, 
$V'_i$ are distinct closed subvarieties of $U'$, $V''_j$ are
distinct closed subvarieties of $U''$, and each $m'_i$ and $m''_j$ is nonzero,
such that $\alpha'|_{U'} = \alpha''|_{U''}$, and glue $\alpha'$ and $\alpha''$ together.
By re indexing the $V'_i$ and the $V''_j$, we can assume that
for $i \le r$, $V'_i \cap U'\cap U'' \ne \emptyset$ and for $i > r$, $V'_i \cap U'\cap U'' = \emptyset$.
Similarly, for 
$j \le s$, $V''_j \cap U'\cap U'' \ne \emptyset$ and for $j > s$, $V''_j \cap U'\cap U'' = \emptyset$.
Note that $V'_i$ can be uniquely recovered from $V'_i \cap U'\cap U''$ if 
$i \le r$, and similarly, $V''_j $
can be uniquely recovered from $V''_j \cap U'\cap U''$ if 
$j \le s$. It follows that $r = s$, and
as closed subvarieties of $X$, and after re-indexing $V'_1,\ldots, V''_r$,
we have $\ov{V_i'} = \ov{V''_i}$ for $i\le r$,
and $m'_i = m''_i$ for $i\le r$ (where $\ov{V}$ denotes the closure of $V$ in $X$).
On the other hand, if $i>r$ so that $V'_i\cap U'' = \emptyset$,
then $V'_i$ is already closed in $X$. Similarly, if $j > r$ then $V''_j$ is already closed in $X$.
Now consider the cycle
$$\alpha = {\textstyle\sum} _{a=1}^r m'_a [\ov{V'_a}] 
+ {\textstyle\sum}_{b=i+1}^r m'_b[V'_b] + {\textstyle\sum}_{c=r+1}^q m''_c[V''_c] \in Z_n(X).$$
It is immediate from its definition that $\alpha|_{U'} = \alpha'$ and
$\alpha|_{U''} = \alpha''$. \hfill $\square$


\bigskip

{\large \bf  4. The sheaf homomorphism from $\K^{\times}$ to $\W$}

\medskip

Let $\Cc^{pure}$ denote the full subcategory of the category $\Cc$ of all algebraic $k$-schemes and
flat $k$-morphisms of constant relative dimensions, for which the objects
are all algebraic $k$-schemes $X$ of pure dimension (which means any two irreducible component of $X$
have the same dimension). We regard $\Cc^{pure}$ as a site under the big Zariski topology.

The sheaf $\Zz_*$ of graded abelian groups on $\Cc^{pure}$, which is the restriction of
the sheaf $\Zz_*$ on $\Cc$ defined in the previous section, has a subsheaf $\W$, which we
call as the sheaf of Weil divisors. This is defined as follows.
For any $X \in \Cc^{pure}$, 
$$\W(X) = Z_{\dim(X) -1}(X)$$
regarded as a graded abelian group concentrated in the degree $\dim(X) -1$.
We also regard $\W$ simply as a sheaf of abelian groups on $\Cc^{pure}$, defined by
the above displayed formula.

For any $X \in \Cc^{pure}$ and any $r \in \K^{\times}(X)$, let
$[r]_X$ denote the Weil divisor of $r$, that is,  
$[r]_X = [\div(r)]$ where $\div(r)$ is the Cartier divisor of $r$. 
By definition,
$$[r]_X = [\div(r)] = {\textstyle \sum}_i \, \ord_{V_i}(r) [V_i]$$
where the sum (which is actually a finite sum as $\ord_{V_i}(r)= 0$ for all but finitely many $V_i$)
is over all closed subvarieties $V_i$ of $X$ of codimension $1$, and
$\ord_{V_i}(r)\in \Z$ is the order of zero (or pole if $< 0$) of $r$ at $V_i$, 
defined as follows. The germ $r_i \in \K_{X,V_i}^{\times}$ of $r$ at the generic point of $V_i$
can be expressed as a ratio $a_i/b_i$ where
$a_i, b_i \in S(\OO_{X,V_i})$ are non zero-divisors, and by definition 
$\ord_{V_i}(r) = \ell(\OO_{X,V_i}/(a_i)) - \ell(\OO_{X,V_i}/(b_i))$ where $\ell(R)$
denotes the length of any artin local ring $R$. This is well-defined, independent
of the choice of $a_i,b_i$, as $\ord_{V_i}$ is an additive function on $S(\OO_{X,V_i})$. 
As the support $|r|$ is closed and of codimension $\ge 1$ in $X$ (see Section 1),
it contains the generic points 
of at most finitely many codimension $1$ subvarieties $V_i$ of $X$, so the above sum
that defines $[r]_X$ is finite.

\bigskip

{\bf Proposition 3.} {\it Let $X \in \Cc^{pure}$ have pure dimension $n$, and
let $X_1,\ldots, X_p$  be the irreducible components of $X$, each regarded as a closed subvariety
of $X$ of dimension $n$. Then the following statements hold.

(1) For any $s \in \S(X)$, we have 
$[s]_X = [ X/(s)] \in  Z_{n-1}(X)$
where $[X/(s)]$ is the fundamental cycle of the 
closed subscheme $X/(s) \subset X$ defined by the ideal sheaf $(s) \subset \OO_X$.

(2) Let $r \in \K^{\times}(X)$, and let 
$r_i  \in \K^{\times}(X_i)$ be its restriction to $X_i$ (see Section 1).
Then we have $[r]_X =
{\textstyle \sum}_i\,  \ell(\OO_{X,X_i}) \, [r_i]_{X_i} \in Rat_{n-1}(X)\subset  Z_{n-1}(X)$. 
} 

\medskip

{\bf Proof.} (1) If $s\in \S(X)$ then 
the support $|s|$ of $s$ is the underlying closed subset $|X/(s)|$
of the closed subscheme $X/(s)$ of $X$. 
Let $V_j$ be the irreducible components of $X/(s)$.
By definition, $[s]_X = {\textstyle \sum}_{V_j} \, \ord_{\OO_{X,V_j}}(s)\,[V_j]
  = {\textstyle \sum}_{V_j} \, \ell(\OO_{X,V_j}/(s)) \,[V_j]$. 
On the other hand, the fundamental cycle $[X/(s)]$ of the scheme $X/(s)$ is given by 
$[X/(s)] = {\textstyle \sum}_{V_j} \, \ell(\OO_{X/(s),V_j}) \,[V_j]$.
As $\OO_{X,V_j}/(s) \cong \OO_{X/(s),V_j}$, the numerical coefficients
of each $V_j$ on the two sides are equal, hence (1) holds.

(2) Now let $r \in \K^{\times}(X)$.
Let $X_1,\ldots, X_p$, where $p\ge 1$ be the irreducible components of $X$,
each regarded as a closed subvariety
of $X$ of dimension $n$. As the support $|r|$ does not contain any height $0$ point
(see Section 1), there are only finitely many closed
subvarieties $V_1,\ldots, V_q$ of dimension $n-1$ of $X$ (where $q\ge 0$) 
that are contained in $|r|$. If $V$ is a closed subvariety of dimension $n-1$ whose generic point
does not lie in $|r|$ then $\ord_V(r) = 0$. Hence
$[r]_X = \sum_{j=1}^q\, \ord_{\OO_{X,V_j}}(r) [V_j]$
(however, it is possible that $\ord_{\OO_{X,V_j}}(r) = 0$
for some $V_j$, even though $V_j \subset |r|$, as the local ring $\OO_{X,V_j}$ may not be regular). 

The cycle $[r_i]_{X_i}$ is by definition the sum 
$[r_i]_{X_i} =    {\textstyle \sum}_{\{ i \, |\, V_j \subset X_i\}}\, \ord_{\OO_{X_i,V_j}}(r_i) \, [V_j]$. 
Hence 
$[r]_X =  {\textstyle \sum}_i\,  \ell(\OO_{X,X_i}) \, [r_i]_{X_i}$ if and only if 
$$ {\textstyle \sum}_j\, \ord_{\OO_{X,V_j}}(r) [V_j] =
{\textstyle \sum}_{\{ (i,j)\, |\, V_j \subset X_i \}}\,  \ell(\OO_{X,X_i})\, \ord_{\OO_{X_i,V_j}}(r_i) \, [V_j]. 
$$
Comparing the numerical coefficients of $[V_j]$ on both sides,
the above holds if and only if for each $j$ we have
$$~~~~~~~~~~~~~~~~~~~~~~\ord_{\OO_{X,V_j}}(r)  =
{\textstyle \sum}_{\{ i\,| \,V_j \subset X_i\}}\, \ell(\OO_{X,X_i})\, \ord_{\OO_{X_i,V_j}}(r_i).
~~~~~~~~~~~~~~~~~~~~\ldots (3)$$
Let $A = \OO_{X,V_j}$, which is noetherian local ring of dimension $1$.
Let $\mm\subset A$ be its maximal ideal and let $\pp_1, \ldots, \pp_k \subset A$ be its
minimal prime ideals (these are defined by the irreducible components $X_i$ which contain $V_j$).
For any $a\in S(A)$ 
the rings $A/(a)$, $A_{\pp_i}$ and 
$A/(\pp_i + (a))$ are artinian. 
By [3] Lemma A.2.7 (applied by taking $M = A$ in the statement of that lemma),
for any $a\in S(A)$ we have
$$\ell(A/(a)) = {\textstyle \sum}_{\pp_i}\, \ell(A_{\pp_i})\, \ell( A/(\pp_i + (a))).$$
In other words, for any $V_j$ and any $a \in S(\OO_{X,V_j})$ we have 
$$~~~~~~~~~~~~~~~~~~~~\ord_{\OO_{X,V_j}}(a)  =
{\textstyle \sum}_{\{ i\, |\, V_j \subset X_i\}}\, \ell(\OO_{X,X_i})\, \ord_{\OO_{X_i,V_j}}(a|_{X_i}).
~~~~~~~~~~~~~~~~~ \ldots (4)$$
If $r'$ denotes the germ of $r$ in $\K_{X,V_j}^{\times}$ 
then $r' = s/t$ for some $s, t \in S(\OO_{X,V_j})$.
Applying (4) to $s$ and $t$ and subtracting the results, we get the
desired equality (3).

As each $r_i$ is a non-zero rational function on the variety $X_i$, by definition
$[r_i]_{X_i}$ belongs to $Rat_{n-1}(X)$, and hence $[r]_X \in Rat_{n-1}(X)$.
\hfill$\square$

\bigskip

{\bf Proposition 4.} {\it The maps $\Phi_X : \K^{\times}(X) \to \W(X)$ defined by 
$\Phi_X(r) = [r]_X$, as $X$ varies over $\Cc^{pure}$, define 
  a homomorphism of sheaves of abelian groups $\Phi : \K^{\times} \to \W$ 
  on $\Cc^{pure}$. The homomorphism $\Phi$ factors via the sub presheaf
  $Rat \subset \W$ defined on $\Cc^{pure}$ by $X \mapsto Rat_{\dim(X) -1}(X) \subset Z_{ \dim(X) -1}(X)$.}

\medskip

{\bf Proof.} 
We just have to show that $\Phi$ commutes with restrictions, that is, 
if $f : X \to Y$ is a flat homomorphism of constant relative
dimension $d$ between pure dimensional schemes $X$ and $Y$ of dimensions $n+d$ and $n$ respectively,
then the following diagram commutes. 
$$\begin{array}{ccc}
\K^{\times}(Y) & \stackrel{f^{\sharp}}{\to} & \K^{\times}(X) \\
  {\scriptstyle \Phi_Y} \da ~~~ && ~~~ \da {\scriptstyle \Phi_X}\\
Z_{n-1}(Y) & \stackrel{f^*}{\to} & Z_{n+d-1}(X) 
\end{array}$$
If $f : X \hra Y$ is an open embedding, then the commutativity of the above diagram
is immediate from the definitions of the objects and the arrows in it.

So now we consider the general case. 
Let $r\in \K^{\times}(Y)$. Then we have two sections 
$$f^*\Phi_Y(r) ,\, \Phi_X(f^{\sharp}r) \in Z_{n+d-1}(X) = \Gamma(X, \Zz_{n+d-1})$$ 
and we have to show that these two are equal, that is, the following equality holds:
$$~~~~~~~~~~~~~~~~~~~~~~~~~~~~~~~~~~~~~~~~~
f^*([r]_Y) = [f^{\sharp}r]_X
~~~~~~~~~~~~~~~~~~~~~~~~~~~~~~~~~~~~~~~~~ \ldots (5)$$ 
As $\Zz_{n+d-1}$ restricts to a sheaf (in particular, a separated presheaf) on the small Zariski site of
$X$, to check the equality (5), by applying the special case of open embeddings,  
it is enough to check that the restrictions of its two sides are equal affine locally on $X$. 

Given any $x\in X$, let $j: V = \spec B \hra Y$ be
an affine open subscheme with $f(x) \in V$. Let $i: U = \spec A \hra X$ be an affine
open subscheme with $x\in U$ and $f(U) \subset V$.
Let $\varphi = f|_U : U\to V$, and let $\varphi^{\sharp} : B \to A$ be the corresponding $k$-algebra
homomorphism. Note that $U$ and $V$ have pure dimensions
$n+d$ and $n$, and $\varphi$ is flat of constant relative dimension $d$.

It follows from the description of meromorphic functions on a noetherian affine scheme
given in Section 1, that $\K^{\times}(V) = Q(B)^{\times}$ and $\K^{\times}(U) = Q(A)^{\times}$,
and the restriction map $\varphi^{\sharp} : \K^{\times}(V) \to \K^{\times}(U)$ is
the ring homomorphism $\varphi^{\sharp} : Q(A)\to Q(B)$ 
induced by 
the flat homomorphism $\varphi^{\sharp} : B \to A$.
The restriction $r|_V = j^{\sharp}r$ of $r \in \K^{\times}(Y)$
to $\K^{\times}(V)= Q(B)^{\times}$ 
can be expressed as a quotient $r|_V = b_1/b_2$ 
where $b_1,b_2 \in S(B)$.
We want to show that the restriction of (5) holds on $U$, that is, 
$$~~~~~~~~~~~~~~~~~~~~~~~~~~\varphi^* ([b_1/b_2]_V) = [\varphi ^{\sharp}(b_1/b_2)]_U  \in Z_{n+d -1}(U)
~~~~~~~~~~~~~~~~~~~~~~~~~~ \ldots (6)$$
Now, $[b_1/b_2]_V = [b_1]_V - [b_2]_V $ and so $\varphi^*([b_1/b_2]_V)
= \varphi ^*([b_1]_V) - \varphi ^*([b_2]_V)$. On the other hand, 
$[\varphi ^{\sharp}(b_1/b_2)]_U = [(\varphi ^{\sharp}b_1)/(\varphi ^{\sharp}b_2)]_U =
[\varphi ^{\sharp}b_1]_U - [\varphi ^{\sharp}b_2]_U$.  
Hence to prove $(6)$, we just have to show 
that for any $b\in S(B)$ the following equality holds:
$$~~~~~~~~~~~~~~~~~~~~~~~~~~~~~~~~~\varphi^* ([b]_V) = [\varphi^{\sharp}b]_U \in Z_{n+d -1}(U)
~~~~~~~~~~~~~~~~~~~~~~~~~~~~~~~~~ \ldots (7)$$
As $b\in S(B) = \S(V)$, by Proposition 3.(1) we have 
$[b]_V = [\spec B/(b)] \in Z_{n-1}(V)$, and so
$\varphi ^*([b]_V) = \varphi ^*[\spec B/(b)]$. 
By definition of flat pullback on cycles, 
$\varphi^*[\spec B/(b)] = [\varphi^{-1} (\spec B/(b))]$, which is the fundamental cycle of
the closed subscheme $\varphi^{-1} (\spec B/(b))$ of $U = \spec A$.
The schematic inverse image $\varphi^{-1} (\spec B/(b))$
is the closed subscheme $\spec A/ (\varphi^{\sharp}b)$ of $\spec A$.
The element $\varphi^{\sharp}b$ lies in $S(A) = \S(U)$ as $\varphi^{\sharp} : B \to A$ is flat.
Hence by Proposition 3.(1) we have 
$[\varphi^{\sharp}b]_U = [\spec A/ (\varphi^{\sharp}b)]$.
Hence we have the sequence of equalities
$$\varphi ^*([b]_V) = \varphi^*[\spec B/(b)] = [\varphi^{-1} (\spec B/(b))]
= [\spec A/ (\varphi^{\sharp}b)] = [\varphi^{\sharp}b]_U$$
in $Z_{n+d -1}(U)$. This proves the equation $(7)$,
and hence completes the proof that $\Phi : \K^{\times} \to \W$ is a homomorphism of sheaves.
The final statement that $\Phi$ factors via the sub
presheaf $Rat \subset \W$ follows from
the equality $[r]_X =  {\textstyle \sum}_i\,  \ell(\OO_{X,X_i}) \, [r_i]_{X_i}$ given by
Proposition 3.(2), because by its definition each $[r_i]_{X_i}$ lies in
$Rat_{\dim(X_i)-1}(X_i) \subset Rat_{\dim(X)-1}(X)$. 
\hfill$\square$

\bigskip

{\bf Remark 5.} We do not need to explicitly
use the morphism $\tilde{b} : \spec B \to \AA^1 = \spec k[t]$
defined by the $k$-homomorphism $k[t] \to B:  t\mapsto b$.
Note that $\tilde{b}^{-1}(P_0) = \spec B/(b)$
where $P_0\in \AA^1$ is the closed point defined by $t=0$.
As $b\in S(B)$, the local criterion for flatness
implies that $\tilde{b}$ is flat in an open neighbourhood of $\spec B/(b)$.
If $B$ is a domain then $\tilde{b}$ is flat if and only if $b$ is transcendental over $k$,
as follows from the nullstellensatz.
The morphism $\tilde{b}$, and implicitly the above flatness condition when $B$ is a domain,  
are implicitly used in the proof 
of the preservation of rational equivalence under flat morphisms in [3].

\bigskip

{\bf Theorem 6.} {\it Let $f : X \to Y$ be a flat morphism of constant relative dimension
  $d$. Then for any $n$, under the homomorphism $f^*: Z_n(Y) \to Z_{n+d}(X)$, the image of
  $Rat_n(Y)$ lies in $Rat_{n+d}(X)$. }

\medskip

{\bf Proof.} By definition, $Rat_n(Y)$ is generated by elements of the form $[r]_V = \Phi_V(r)$ 
where $V\subset Y$ is a closed subvariety of dimension $n+1$ and $r\in \K(V)^{\times}$.
Hence we just have to show that $f^*\Phi_V(r)$ lies in $Rat_{n+d}(X)$.
Let $W = f^{-1}(V) \subset X$ be the schematic inverse image of $V$.  
Hence $g = f|_W : W \to V$ is again flat and of constant relative
dimension $d$. As $V$ is a variety hence of pure dimension $n$, this implies that 
$W$ is a scheme of pure dimension $n+d+1$, and so both $V,\,W\, \in \Cc^{pure}$,
and hence $g$ is an arrow in $\Cc^{pure}$, which is a full subcategory of $\Cc$. 
Hence by Proposition 4, 
$$g^*\Phi_V(r) =  \Phi_W(g^{\sharp} r)$$
where $g^{\sharp} : \K^{\times}(V) \to \K^{\times}(W)$ is the restriction homomorphism 
under $g$ for the sheaf $\K^{\times}$ over $\Cc^{pure}$.
By Proposition 3.(2) applied to the scheme $W$ and to $g^{\sharp} r \in \K^{\times}(W)$,
the cycle $\Phi_W(g^{\sharp} r) = [g^{\sharp} r]_W$ lies in the subgroup 
$Rat_{n+d}(W) \subset Rat_{n+d}(X)$. Hence the theorem.
\hfill $\square$

\bigskip

{\bf Remark 7.} The following rough analogy partly inspired our sheaf theoretic approach.
The sub presheaf $Rat \subset \W$
on the site $\Cc^{pure}$,  defined by $X \mapsto Rat_{\dim(X)-1}(X)$,
is not a sheaf because a locally principal Weil divisor 
need not be rationally equivalent to $0$ globally. This means that the presheaf $X\mapsto A_*(X)$
on $\Cc$ is not separated, and hence not a sheaf. 
The situation is similar to that for differential forms
on manifolds. Let $\M$ be the site formed by all $C^{\infty}$ manifolds
and $C^{\infty}$ maps, equipped with the big Euclidean topology. For any integer $n$,
differential $n$-forms form a sheaf $\Omega^n$ on the site $\M$, and 
closed differential $n$-forms form a subsheaf $\Zz^n$ of $\Omega^n$.
The exact $n$-forms only form a sub presheaf $\B^n \subset \Zz^n$, which
however is not a sheaf as a closed form is locally exact but not necessarily globally exact.
The de Rham cohomology $H_{dR}^*$ forms a presheaf which is not separated, so not a sheaf.
Just as we can define $C^{\infty}$-pullbacks on $H_{dR}^*$ using the sheaf homomorphism
$d : \Omega^{n-1} \to \Zz^n$, one should be able to define flat pullbacks on $A_*$ using the sheaf 
homomorphism $\Phi : \K^{\times} \to \W$ on $\Cc^{pure}$, with $Rat$ playing the role of $\B^n$.
This appears to work.

\medskip

{\it Postscript.} It was brought to my attention that our 
sheaf theoretic approach to Chow groups has certain
similarities with the viewpoint of Rost [6], which
can also be used to derive a parallel proof of the preservation of rational equivalence
under flat pullbacks.


\bigskip


{\bf\large  References}

[1] C. Chevalley, Les classes d'\'equivalence rationelle, I, II. S\'eminaire C. Chevalley,
2$^e$ ann\'ee, {\it Anneaux de Chow et Applications}, Secr. Math. Paris, 1958.

[2] A. Grothendieck, A. and J. Dieudonne, {\it \'El\'ements de G\'eom\'etrie Alg\'ebrique
- IV$_4$}, Publ. Math. IHES 32 (1967). 

[3] W. Fulton, {\it Intersection Theory}, Springer Verlag, 1984.

[4] S. Kleiman, Misconceptions about $K_X$, Enseign. Math. 25 (1979) 203-206.

[5] D. Mumford, {\it Lectures on Curves on an Algebraic Surface}, Princeton University Press, 1966.

[6] M. Rost, {\it Chow groups with coefficients}, Doc. Math. J. DMV 1 (1996) 319-393.


\end{document}